\documentclass[12pt]{amsart}

\usepackage{amsfonts,amsmath,latexsym,amssymb,verbatim,amsbsy}
\usepackage{amsthm}

\usepackage{pstricks}

\theoremstyle{plain}
\newtheorem{THEOREM}{Theorem}[section]

\newtheorem{theorem}[THEOREM]{Theorem}

\newtheorem{lemma}[THEOREM]{Lemma}

\theoremstyle{definition}

\theoremstyle{remark}

\newcommand{\thm}[1]{Theorem~\ref{#1}}

\def \a {\alpha}
\def \b {\beta}

\def \f {\varphi}

\def \k {\kappa}

\def \n {\nabla}

\def \O {\Omega}

\def \LH {\mathcal{LH}}

\def \< {\langle}
\def \> {\rangle}
\def \p {\partial}
\def \ra {\rightarrow}

\newcommand{\be}{\begin{equation}}
\newcommand{\ee}{\end{equation}}
\newcommand{\uh}{{u^\mathrm{h}}}
\newcommand{\ul}{{u^\mathrm{l}}}

\begin{document}

\title[Energy equality for 3D Navier-Stokes]{On the energy equality for weak solutions of the 3D Navier-Stokes equations}

\author{ A. Cheskidov}
\address[A. Cheskidov]
         {Department of Mathematics\\
          University of Michigan\\
          Ann Arbor, MI 48109}
\email{acheskid@umich.edu}

\author{S. Friedlander}
\address[S. Friedlander and R. Shvydkoy]
{Department of Mathematics, Stat. and Comp. Sci.\\
        University of Illinois\\
        Chicago, IL 60607}
\email{susan@math.northwestern.edu}

\author{R. Shvydkoy}
\email{shvydkoy@math.uic.edu}

\date{\today}

\begin{abstract}
We prove that the energy equality holds for weak solutions of the 3D
Navier-Stokes equations in the functional class $L^3([0,T);V^{5/6})$, where
$V^{5/6}$ is the domain of the fractional power of the Stokes operator $A^{5/12}$.
\end{abstract}



\maketitle

\section{Introduction}
We consider the 3D incompressible Navier-Stokes equations (NSE)
\begin{equation} \label{NSE1}
\left\{
\begin{aligned}
&\p_t u - \nu \Delta u + (u \cdot \nabla)u + \nabla p = f,\\
&\nabla \cdot u =0 \\
&u(x,t) = 0 \text{ for } x \in \p \Omega,
\end{aligned}
\right.
\end{equation}
on an open bounded domain  $\Omega \subset \mathrm{R}^3$ of class $C^2$.
Here $u(x,t)$, the velocity, and $p(x,t)$, the pressure, are unknowns; $f(x,t)$ is
a given driving force which we assume to belong to $L^1([0,T);L^2(\Omega))$, and $\nu>0$ is the kinematic  viscosity coefficient of the fluid.

By the classical result of Leray \cite{Leray} and Hopf \cite{Hopf},
for every divergence-free initial data $u_0 \in L^2(\Omega)$ and
$T>0$ there exists a weak solution to the system \eqref{NSE1} in the
class
\begin{equation}
\LH = L^\infty_{\rm loc}([0,T);L^2(\Omega)) \cap L^2_{\rm loc}([0,T);H^1(\Omega)),
\end{equation}
satisfying
\begin{multline}\label{weakform}
 \int_0^T \left\{ - (u,\p_t \f) + \nu (\n u,\n \f) + (u\cdot \n u, \f) \right\} dt\\ = (u_0,\f(0)) + \int_0^T (f,\f)dt,
\end{multline}
for all test functions $\f \in C_0^\infty([0,T)\times \Omega)$ with $\n \cdot \f = 0$.
In addition one can find a weak solution satisfying the strong energy inequality, i.e.
\begin{equation}\label{SEI}
|u(t)|_2^2 + 2\nu \int_{t_0}^t |\n u(s)|_2^2 ds \leq |u(t_0)|_2^2 + \int_{t_0}^t (f(s) \cdot u(s)) ds,
\end{equation}
for all $t \in [0,T)$ and almost all $t_0 \in [0,t)$ including $t_0
= 0$. With an additional correction of $u$ on a subset of $[0,T)$ of
measure zero one can ensure that $u$ is weakly continuous with
values in $L^2(\Omega)$ (see Serrin \cite{Ser63}), and
\eqref{weakform} holds in a stronger form,
\begin{multline}\label{weakform2}
 (u(t),\f(t)) + \int_0^t \left\{ - (u,\p_s \f) + \nu (\n u,\n \f) + (u\cdot \n u, \f) \right\} ds\\ = (u_0,\f(0)) + \int_0^t (f,\f)ds,
\end{multline}
for all $t\in[0,T)$ and $\f$ as before.

Solutions satisfying all the properties
listed above are commonly called Leray-Hopf solutions. 
It is a famous open problem to show that
a given Leray-Hopf solution to the Navier-Stokes system actually satisfies the energy equality
\begin{equation}\label{EE}
|u(t)|_2^2 + 2\nu \int_{t_0}^t |\n u(s)|_2^2 ds = |u(t_0)|_2^2 + \int_{t_0}^t (f(s) \cdot u(s)) ds,
\end{equation}
for all $0\leq t_0 \leq t < T$ (or equivalently all $t$ and $t_0
= 0$). Continuing interest to this problem is motivated by the fact
that the failure of \eqref{EE} opens a possibility for an energy
sink other than the natural viscous dissipation. Such a property of
real fluids is not expected to exist physically. It has therefore
been the subject of many works in the past and in recent years to
find sufficient conditions for the energy equality to hold in
bounded and unbounded domains. Below we mention the results that are
most relevant to ours.

Serrin \cite{Ser63} showed \eqref{EE} under the assumption that 
\begin{equation}
 u\in \LH \cap L^r([0,T);L^s), \text{ for } s \geq 3,\ \frac{3}{s}+\frac{2}{r} = 1.
\end{equation}
Lions and Ladyzhenskaya improved the scaling in Serrin's condition to $u
\in L^4([0,T);L^4)$ (see \cite{LSU,Lions}). Shinbrot in \cite{Shinbrot}
proved \eqref{EE} under the extrapolated version of the
Lions-Ladyzhenskaya condition, namely $2/s+2/r \leq 1$ for $s \geq
4$. One can show that Shinbrot's condition in the functional class
$\LH$ is weaker that that of Serrin (see \cite{Galdi}). A
dimensionally the same, yet formally weaker, condition was recently
found by Kukavica \cite{Kuk06}. He shows that \eqref{EE} holds if
the pressure $p$ is locally $L^2$-summable in space-time. Since $p$
is given by a Calderon-Zygmund operator applied to $u\otimes u$,
Kukavica's condition is implied by that of Lions and Ladyzhenskaya.

In this present paper we prove \eqref{EE} in a functional class 
with a better scaling than that of Shinbrot. Let us denote
\begin{equation}\label{H}
H=\{u\in L^2(\O): \nabla\cdot u =0, u\cdot n|_{\partial\O}=0\},
\end{equation}
and let $\mathbb{P} : L^2(\Omega) \to H$ be the $L^2$-orthogonal
projection, referred to as the Helmholtz-Leray projector. Let $A$ be
the Stokes operator defined by \be\label{A} Au = -\mathbb{P}\Delta
u. \ee The Stokes operator is a self-adjoint positive sectorial
operator with a compact inverse. We denote $V^{s} =
\mathcal{D}(A^{s/2})$, $s>0$, the domain of the fractional power of
$A$ (see the next section for details). Roughly, $V^{s}$ corresponds
to the Sobolev space $H^s$, and is in fact equal to $H^{s}$ under
periodic boundary conditions. In particular one can give the
following description of $V=V^1$ (see \cite{CF}): 
\be\label{V}
V=\{u\in H^1(\O): \nabla\cdot u =0, u|_{\partial\O}=0\}. 
\ee 
We prove the following result.
\begin{theorem}\label{T:main}
Every weak solution $u(t)$ of \eqref{NSE1} satisfying \eqref{weakform2} on $[0,T)$ with $u \in
\LH \cap L^3([0,T); V^{5/6})$ verifies the energy equality
\eqref{EE}.
\end{theorem}
The dimensional $L^rL^s$-analogue of the functional class
$L^3([0,T); V^{5/6})$ is $L^3L^{9/2}$, which exhibits the scaling
\begin{equation}\label{new}
2/s+2/r = 10/9.
\end{equation}
It would be an interesting next step to try to prove the energy equality under this new
scaling in the $L^rL^s$-class itself.

Finally we note that under periodic boundary conditions the result
of \thm{T:main} follows from  -- although not explicitly stated in
-- the recent work by the authors on Onsager's conjecture for the
Euler equations \cite{ccfs}. So, the novelty of this present paper
consists in the extension of \thm{T:main} to the case of the
Dirichlet boundary condition in a smooth domain.

\section{Preliminaries}

First, let us introduce some notations and functional setting.
Denote by $(\cdot,\cdot)$ and $|\cdot|$ the $L^2(\Omega)$-inner product and the
corresponding $L^2(\Omega)$-norm.
Let $H$ and $V$ be as above \eqref{H}, \eqref{V}, and $V'$ denote the dual of $V$.
We endow $V$ with the norm $\|u\| = |\n u|$.
There exists an orthonormal basis of eigenvectors $\{w_n\}$ in $H$, and a sequence of positive
eigenvalues $\{\lambda_n\}$, such that
\be
Aw_n = \lambda_n w_n, \qquad w_n \in \mathcal{D}(A),
\ee
and
\be
0<\lambda_1 \leq \lambda_2 \leq \dots \leq \lambda_n \leq \dots,
\qquad \lim_{n \to \infty} \lambda_n = \infty.
\ee
Any $u\in H$ can then be written as
\be
u = \sum_{n=1}^\infty (u,w_n) w_n.
\ee
Henceforth we will use the notation $u_n=(u,w_n)$.
For $s>0$ we define the operator $A^s$ by
\be
A^{s}u = \sum_{n=1}^\infty \lambda_n^s u_n w_n,
\ee
and the space
\be
V^{s}=\{u \in H: u=\sum_{n=1}^\infty u_n w_n,\  \|u\|_s^2 = \sum_{n=1}^\infty \lambda_n^{s} |u_n|^2<\infty \}.
\ee
Thus, $V^{s} = \mathcal{D}(A^{s/2})$.

Now denote $B(u,v):=\mathbb{P}(u \cdot \nabla v)\in V'$ for
$u, v \in V$. So, we can rewrite (\ref{NSE1}) as the following differential
equation in $V'$:
\begin{equation} \label{NSE}
\p_t u + \nu A u +B(u,u) = g,\\
\end{equation}
where $u$ is a $V$-valued function of time and $g = \mathbb{P} f$.
Finally, we denote $b(u,v,w) = \langle B(u,v),w \rangle$.
This trilinear form is anti-symmetric:
\[
b(u,v,w)=-b(u,w,v), \qquad u,v,w \in V,
\]
in particular, $b(u,v,v)=0$ for all $u,v \in V$.

\section{The proof of \thm{T:main}}

Define
\be
P_\k u= \sum_{n:\lambda_n \leq \kappa} u_n w_n, \qquad u \in H.
\ee
Let $u\in V^{\beta}$ and denote $\ul=P_\k u$, $\uh=u-\ul$.
Observe the following inequalities:
\begin{equation}\label{ba}
\begin{aligned}
\|\ul\|_\b &\leq \k^{\b - \a} \|\ul\|_\a  \\
\|\uh\|_\a &\leq \k^{\a - \b} \|\uh\|_\b,
\end{aligned}
\end{equation}
whenever $\b > \a$.

\begin{lemma} \label{tr}
Let $u(t)$ be a weak solution of \eqref{NSE1} on $[0,T)$. Then
\begin{multline}
    |u(t)|^2 +  2\nu\int_{t_0}^t \|u\|^2 \, ds \\= |u(t_0)|^2 + 2\int_{t_0}^t (g, u) \, ds +
      2\lim_{\k \to \infty} \int_{t_0}^t b(u,\ul,u) \, ds,
\end{multline}
for all $0 \leq t_0 \leq t < T$.
\end{lemma}
\begin{proof}
One can see directly from \eqref{weakform2} that $\ul \in
C_w([0,T); V)$ and $\p_t \ul \in L^2_{loc}([0,T);V)$. Thus, using
$\ul$ as a test function in \eqref{weakform2} (allowed by a standard
approximation argument)  we obtain
\begin{multline}
    |\ul(t)|^2 - |\ul(t_0)|^2 + 2\nu\int_{t_0}^t \|\ul\|^2 \, ds
    - 2\int_{t_0}^t (g, \ul) \, ds \\ =2\int_{t_0}^t b(u,\ul,u) \, ds.
\end{multline}
>From this we see that the limit of the right hand side exists as $\k
\to \infty$, which completes the proof of the lemma.
\end{proof}

Let $\ul$ and $\uh$ be defined as before.
In view of Lemma \ref{tr}, it suffices to show that
\be
\lim_{\k \to \infty} \int_0^T |b(u,\ul,u)| \, ds = 0.
\ee
Indeed, let us write
$$
b(u,\ul,u) = b(\uh,\ul,\uh) + b(\ul,\ul,\uh) +
b(\uh,\ul,\ul)+b(\ul,\ul,\ul).
$$
The last two terms vanish, so it suffices to estimate only the first
two. We use the standard estimate found, for example, in \cite{CF}:
\begin{equation}\label{}
    |b(u,v,w)| \leq \|u\|_{s_1} \| v\|_{s_2+1} \|w\|_{s_3}
\end{equation}
where $s_1+s_2+s_3 \geq 3/2$. To estimate the first term let us set
$s_1 = s_2 = s_3 = 1/2$, then
$$
|b(\uh,\ul,\uh)| \leq \|\uh\|^2_{1/2} \|\ul\|_{3/2},
$$
and by \eqref{ba} we have
\begin{align}
\|\uh\|_{1/2} &\leq  \k^{-1/3}\|\uh\|_{5/6} \\
\|\ul\|_{3/2} & \leq  \k^{2/3} \|\ul\|_{5/6}.
\end{align}
So,
$$
|b(\uh,\ul,\uh)| \leq  \|\uh\|_{5/6}^2\|\ul\|_{5/6},
$$
which tends to zero a.e. in $t$ as $\k \ra \infty$. Since in
addition,
$$
|b(\uh,\ul,\uh)| \leq  \|u\|^3_{5/6}
$$
for all $t$, by the Dominated Convergence Theorem,
$$
|b(\uh,\ul,\uh)| \to 0, \qquad \text{as} \qquad \k \to \infty,
$$
in $L^1([0,T))$. As to the second term, similar estimates with $s_1
= 5/6$, $s_2 = 0$, $s_3 = 2/3$, yield
$$
|b(\ul,\ul,\uh)| \leq  \|\ul\|^2_{5/6}\|\uh\|_{5/6},
$$
which also tends to zero in $L^1([0,T))$ as $\k \ra \infty$. Hence
we obtain the desired energy equality.

\section*{acknowledgment}
The work of AC was partially supported by NSF PHY grant 0555324, the
work of SF by NSF DMS grant 0503768, and the work of RS by NSF DMS
grant 0604050. The authors thank V. \v{S}ver\'{a}k for useful
discussions.


\end{document}